\documentclass[12pt]{amsart}

\usepackage{PascalianStyle}
\title{Polynomials Arising from Sorted Binomial Coefficients}

\author[Owen Levens]{Owen Levens$^*$}
\address{Department of Mathematical Sciences, DePaul University, Chicago, IL, USA}
\email{olevens@depaul.edu}
\thanks{$^*$Research partially supported by the DePaul University Undergraduate Summer Research Program.}

\begin{document}

\maketitle
\begin{abstract}
    The triangle of sorted binomial coefficients $\B{n}{k}=\binom{n}{\lfloor (n-k)/2\rfloor}$ for $0\leq k\leq n$ has appeared several times in recent combinatorial works but has evaded dedicated study. Here we refer to $\B{n}{k}$ as the \textit{Pascalian numbers} and unify the various perspectives of $\B{n}{k}$. We then view each row of the $\B{n}{k}$ triangle as the coefficients of the \textit{$n$th Pascalian polynomial}, which we denote $P_n(z)$. We derive recursions, formulae, bounds on $P_n(z)$'s roots in $\CC$ and characterize the asymptotics of these roots. We show the roots of $P_n(z)$ converge uniformly to a curve $\partial\Gamma\in\CC$ and asymptotically fill the curve densely. We conclude with a discussion of the reducibility and Galois groups of $P_n(z)$. Our work has natural connections to the truncated binomial polynomials, asymptotic analysis, and well known integer families.
\end{abstract}

\section{Introduction}
Various recent combinatorial interpretations \cite{ActonPetersenShirmanTenner2025,levens2025} have naturally emerged for the triangle of numbers given by $\binom{n}{\lfloor (n-k)/2\rfloor}$ \cite[A061554]{OEIS}, which has sporadically appeared previously in the Riordan-array literature \cite{hennessy2011,sunma2014,cigler2015,guprodinger2021}. We unify and expand certain perspectives on these numbers then initiate the analytic, asymptotic, and algebraic studies of the associated polynomials. Particularly, we let $\B{n}{k}=\binom{n}{\lfloor (n-k)/2\rfloor}$ for $0\leq k\leq n$, call the $\B{n}{k}$ triangle  the \textit{Pascalian numbers}, and let the $n$th \textit{Pascalian polynomial} be

$$P_n(z)=\sum_{k=0}^n\B{n}{k}z^{n-k}.$$ 
Our results relate to ideas beyond $P_n(z)$'s study. For instance, evaluating our expression for $P_n(z)$'s generating function at specific values $z_0$ yields the combinatorial generating functions encoding the sequence $P_n(z_0)$. Similarly, our method of proving the prime factorization of $P_{2m+1}(z)$ applies to a far wider class of polynomials. 

Section $\ref{sec: combinatorial foundations}$ formally introduces $\B{n}{k}$ from the seperate combinatorial perspectives of \cite{ActonPetersenShirmanTenner2025} and \cite{levens2025} then derives the formula $\B{n}{k}=\binom{n}{\lfloor(n-k)/2\rfloor}$--ultimately leading to a new correspondence between two row standard domino tableaux and certain lattice walks which preserves multiple statistics. The Pascalian numbers are then realized on binary strings. In Section $\ref{sec: PP}$, we turn to the Pascalian polynomials themselves. We derive several recursions, formulae, and the generating function for $P_n(z)$.

We turn our attention to the roots of the Pascalian polynomials for the next two sections. Section $\ref{sec: roots of PP}$ concerns non-asymptotic questions of $P_n(z)$'s roots, such as the number of real and purely imaginary roots of $P_n(z)$ and bounds in $\CC$. We show $P_n(z)$ has at most one real root (called the \textit{trivial root}), a unique pair of purely imaginary roots if and only if $n\equiv 3\mod 4$, and bound the roots within the optimal annulus in $\CC$.
\vspace{.1in}

\noindent\textbf{Theorem~\ref{thm: PP roots annuli}.} \textit{The nontrivial roots of $P_n(z)$ lie within the annulus $\sqrt{2}-1<|z|<1$.}

\vspace{.1in}

\noindent\textbf{Lemma \ref{lem: unique real root}.} \textit{The trivial root at $-1$ is the unique real root of $P_n(z)$ for odd $n$. For even $n$ and real $z$, $P_n(z)$ is strictly positive.}

\vspace{.1in}

We use these results, with a recursion of Section \ref{sec: PP}, to show that $P_n(z)$ and $P_{n-2}(z)$ may only share the trivial root in Section $\ref{sec: roots of PP}$ before passing to asymptotic considerations. In Section $\ref{sec: limit curve of roots}$, we introduce a family of curves to help describe the root asymptotics of $P_n(z)$. Particularly, we set $n\geq 2$,
$$K_n=\frac{n^2-1}{2n^2},$$
$$\Gamma_n = \left\{ z\in\CC:\ \frac{|z|}{\sqrt[n]{|1+z|}|1+z^2|}\leq K_n,\ |z|\le 1\right\},$$
and $\partial\Gamma_n$ to the boundary of $\Gamma_n$. With $\Gamma$ equal to the limit of the $\Gamma_n$ and $\partial\Gamma$ to its boundary, we prove the following theorems.

\vspace{.1in}
\noindent\textbf{Theorem \ref{thm: Pn(z) has no roots in D_n}.} \textit{$P_n(z)$ has no roots in $\Gamma_n$ for all $n\geq 2$.}

\vspace{.1in}
\noindent\textbf{Theorem \ref{thm: PP roots converge to C_infty}.} \textit{The roots of $P_n(z)$ converge uniformly to $\partial\Gamma$ and fill the curve densely.}
\vspace{.1in}

Section $\ref{sec: irreducibility and galois groups}$ discusses on the reducibility and Galois groups of $P_n(z)$. For $n=2m+1$, we give $P_{n}(z)$'s factorization into irreducibles and embed its Galois group inside $S_{m}^B$, the $m$th hyperoctahedral group. Several questions concerning the reducibility and Galois groups of $P_n(z)$ are left open as conjectures, particularly for even $n$. Natural connections emerge to active work on the truncated binomial polynomials. We close with various directions for further study and open questions.

\section{Combinatorial Foundations}\label{sec: combinatorial foundations}

The combinatorial construction of $\B{n}{k}$ from \cite{levens2025} will be our starting point. From there, we derive key facts and present a bijection with the objects given in \cite{ActonPetersenShirmanTenner2025} which unifies multiple notions from the independent contexts. 

\subsection{Standard Domino Tableaux}\label{subsec: SDT}

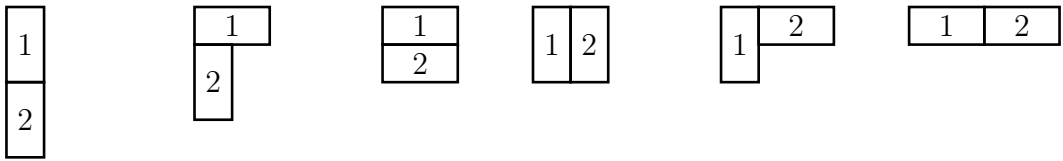
\begin{figure}[h]
\centering
\begin{tikzpicture}[line width=1pt, scale=0.5]

\draw (-12,0) rectangle (-11,2);  \node at (-11.5,1) {1};
\draw (-12,-2) rectangle (-11,0); \node at (-11.5,-1) {2};

\draw (-7,1) rectangle (-5,2);    \node at (-6,1.5) {1};
\draw (-7,-1) rectangle (-6,1);   \node at (-6.5,0) {2};

\draw (-2,1) rectangle (0,2);     \node at (-1,1.5) {1};
\draw (-2,0) rectangle (0,1);     \node at (-1,0.5) {2};

\draw (2,0) rectangle (3,2);      \node at (2.5,1) {1};
\draw (3,0) rectangle (4,2);      \node at (3.5,1) {2};

\draw (7,0) rectangle (8,2);      \node at (7.5,1) {1};
\draw (8,1) rectangle (10,2);     \node at (9,1.5) {2};

\draw (12,1) rectangle (14,2);    \node at (13,1.5) {1};
\draw (14,1) rectangle (16,2);    \node at (15,1.5) {2};
\end{tikzpicture}
    \caption{The $6$ standard domino tableaux with $2$ dominos. From left to right, they have shape $(1,1,1,1), (2,1,1), (2,2), (2,2),( 3,1),$ and $(4)$.}
    \label{fig: 2 domino SDT}
\end{figure}

A weekly decreasing list of integers $\lambda=(\lambda_1,\lambda_2,\dots,\lambda_k)$ is said to partition $n$ if $\sum_{j=0}^k \lambda_j=n$, in which case we write $\lambda\vdash n$. We construct a \textit{domino diagram of shape $\lambda\vdash 2n$} by placing $n$ non overlapping $1\times 2$ dominos in the plane with $\lambda_i$ boxes filling row $i$. Since $\lambda$ is an integer partition, all domino tableaux are upper left justified. We call a domino diagram of $n$ dominos a \textit{Standard Domino Tableau} (SDT) if the dominos are uniquely labeled $1,2,\dots n$ and values appear in increasing order across rows and down columns (see Figure $\ref{fig: 2 domino SDT}$). We call the domino labeled $x$ the $x$ domino and let SDT$_n$ denote the set of SDT with $n$ dominos.

The Garfinkle–Barbasch–Vogan (GBV) correspondence is a bijection between pairs in SDT$_n$ of the same shape and elements of the $n$th hyperoctahedral group $S_n^B$--often represented as the group of size $n$ signed permutations (see Figure $\ref{fig:gbv-tableaux-quarter}$). 
Refer to Leeuwen \cite{vanleeuwen} for a combinatorial treatment under an RSK perspective or Barbasch and Vogan \cite{barbaschvogan} or Garfinkle \cite {garfinkle} for its Lie theoretic origins.
\begin{figure}[h]
\centering
\[
w = -3, 1, 2
\quad \overset{\text{GBV}}{\longleftrightarrow} \quad
\left(
\begin{tikzpicture}[baseline=-0.5ex, x=0.5cm, y=0.5cm, line width=0.6pt]
  \draw (0,0) rectangle (2,1); \node at (1,0.5) {1};
  \draw (2,0) rectangle (4,1); \node at (3,0.5) {2};
  \draw (0,-1) rectangle (2,0); \node at (1,-0.5) {3};
\end{tikzpicture}
\; , \;
\begin{tikzpicture}[baseline=-0.5ex, x=0.5cm, y=0.5cm, line width=0.6pt]
  \draw (0,-1) rectangle (1,1); \node at (0.5,0) {1};
  \draw (1,-1) rectangle (2,1); \node at (1.5,0) {2};
  \draw (2,0) rectangle (4,1); \node at (3,0.5) {3};
\end{tikzpicture}
\right)
\]
\caption{$-3,1,2\in S_3^B$ represented as a pair of SDT of shape $(4,2)$}
\label{fig:gbv-tableaux-quarter}
\end{figure}

In \cite{levens2025}, the GBV correspondence proved useful in the study of global patterns of signed permutations by permitting a signed analogue of Greene's Theorem (see \cite[Ch. 7]{Stanley1999}). Global patterns of signed permutations formally appeared in \cite{lewistenner} and were further studied in \cite{levens2025}--but had previously emerged  under different language (see Egge's work in \cite{egge2007,egge2010}, for instance). Like the authors of \cite{levens2025}, we naturally view each size $n$ signed permutation $w(1),w(2),\dots w(n)$ as size a $2n$ strings $w(-n),w(-n+1),\dots w(-1),w(1),\dots w(n)$ with $w(-i)=-w(i)$. Any size $k$ increasing/decreasing subsequence of this size $2n$ string is called a global $12\dots k$ pattern or a global increasing/decreasing subsequence--a special case of the more general notion of global patterns. They showed that, if $w\in S_n^B$ corresponds to two SDT of shape $\lambda=(\lambda_1,\lambda_2,\dots \lambda_k)\vdash 2n$ under GBV, then $k$ is the length of the longest decreasing global subsequence and $\lambda_1$ is the length of the longest increasing global subsequence\footnote{There are various notations for these objects, all equivalent up to rotation and reflection. With another choice of notation $\lambda_i$ and $k$ may need to be flipped.}. 
For instance, $-3,1,2\in S_3^B$ has largest increasing global subsequence $-2,-1,1,2$ (size $4$) and several largest decreasing global subsequences of size $2$, such as $3,-3$. Then $w$ corresponds to a pair of SDT of shape $(4,2)$, as demonstrated in Figure $\ref{fig:gbv-tableaux-quarter}$.

Permutation pattern avoidance is a major and active field of Combinatorics, with significant contributions from Vincent Vatter \cite{vatter-classes}, Sergey Kitaev (see \cite{kitaev-ppw} for a monograph), and Bridget Tenner \cite{ActonPetersenShirmanTenner2025,levens2025,lewistenner}, who maintains the Database of Permutation Pattern Avoidance \cite{dppa}. The authors of \cite{levens2025} used the GBV correspondence to provide $2$ novel proofs that there are $\binom{2n}{n}$ size $n$ signed permutation with no global $321$ pattern--with the same enumeration for those with no global $123$ pattern. This enumeration has seen several proofs in many distinct contexts. The authors of \cite{levens2025} showed the desired signed permutations are Stembridge's top fully commutative elements, which makes his count of these in \cite{stembridge} as $\binom{2n}{n}$ the first proof. Egge later gave two, more direct proofs in \cite{egge2007,egge2010}, the later of which implicitly used \cite{levens2025}'s signed Greene's theorem result in a distinct way from both proofs of \cite{levens2025}.
We outline \cite{levens2025}'s first argument below. 

Let $B_n$ denote the subset of SDT$_n$ fitting in $2$ rows, since pairs of these represent the signed permutations with no global $321$ pattern (see the last four SDT in Figure \ref{fig: 2 domino SDT} for $B_2$). Particularly, the number of pairs in $B_n$ of shape $(n+k,n-k)$ for $k=0,1,\dots,n$ gives the number of signed permutations globally avoiding $321$. Once we know how many SDT in $B_n$ have shape $(n+k,n-k)$ for each $k$, we sum their squares to count pairs. Note that $0\leq k\leq n$ to respect upper justification. 
\begin{definition}\label{def: pascalian numbers}
    The set of tableaux in $B_n$ with shape $(n+k,n-k)$ is denoted $B(n,k)$. We call the sizes of $B(n,k)$ the \textit{Pascalian numbers} and denote $|B(n,k)|$ as $\B{n}{k}$.
\end{definition}
\begin{figure}[h]
    \centering
\begin{gather*}
    1\\
    1 \hspace{.2 in} 1\\
    2 \hspace{.2 in} 1 \hspace{.2 in} 1\\
    3 \hspace{.2 in} 3 \hspace{.2 in} 1 \hspace{.2 in} 1\\
    6 \hspace{.2 in} 4 \hspace{.2 in} 4 \hspace{.2 in} 1 \hspace{.2 in} 1\\
    10 \hspace{.2 in} 10 \hspace{.2 in} 5 \hspace{.2 in} 5 \hspace{.2 in} 1 \hspace{.2 in} 1 \\
\end{gather*}
    \caption{Pascalian numbers represented in a triangular array for $0\leq n\leq5$.}
    \label{fig: SDT with 2 dominos}
\end{figure}
Subtracting $\lambda_1-\lambda_2$ shows the row lengths may only differ by a multiple of $2$. We see this combinatorially by appending an $n+1$ domino to elements of $B_n$: either a single row grows by $2$ or both extend by $1$, preserving an even difference regardless. As Figure $\ref{fig: SDT with 2 dominos}$ demonstrates for small values, \cite[Thm. 5.1]{levens2025} showed
the triangular array made by $\B{n}{k}$ agrees with pascal's triangle after sorting rows into decreasing order \cite[A061554]{OEIS}.

\begin{theorem}\label{thm: 2 col SDT satisfy pascalian recursion}
There are $\binom{n}{\lfloor\frac{n-k}{2}\rfloor}$ SDT of shape $(n+k,n-k)$, so $\B{n}{k}=\binom{n}{\lfloor\frac{n-k}{2}\rfloor}$. Particularly, the Pascalian numbers satisfy the recursion
 $$\B{n}{k}=
 \begin{cases}
\B{n-1}{k-1}+\B{n-1}{k+1} \text{ for }0<k< n-1\\
\B{n-1}{0}+\B{n-1}{1} \text{ for }k=0
\end{cases}$$
with $\B{n}{n}=\B{n}{n-1}=1$ for each $n$.
\end{theorem}

Then the number of size $n$ signed permutations with no size $3$ increasing global subsequence is the following 
$$\sum_{k=0}^n\B{n}{k}^2=\sum_{k=0}^n\binom{n}{k}^2=\binom{2n}{n}$$
where the first equality is a reordering terms. Additionally, $\B{n}{0}+\B{n}{1}\dots+\B{n}{n}=2^n$. To see this combinatorially, we consider how elements of $B(n,k)$ extend into larger tableaux by concatenating an $n+1$ domino. If $k=0$, we may either place a vertical $n+1$ domino or a horizontal $n+1$ domino in row $1$. When $0<k$, we extend each element of $B(n,k)$ by appending a horizontal $n+1$ domino in either row. Notice that, since $\lambda_1$ and $\lambda_2$ differ by a multiple of $2$, we may always place a horizontal $n+1$ domino in row $2$ without violating upper justification. Thus $|B_n|=2|B_{n-1}|\dots =2^{n-1}|B_1|=2^{n}$.

With $2^n$ objects, we might wonder if an alterative statistic on $B_n$ gives a partition of $B_n$ with the more natural distribution $\binom{n}{k}$. To give this statistic on $B_n$, we have the following map: for $T\in B_n$, let $S_T\subseteq[n]$ be the set of labels of dominos which appear horizontally in row $1$ of $T$. 
\begin{theorem}\label{thm: ST determines T}
    The map $T\mapsto S_T$ is a bijection between elements of $B_n$ and subsets of $[n]$, so elements of $B_n$ are uniquely determined by the dominos appearing horizontally in row $1$. Further, there are $\binom{n}{k}$ elements of $B_n$ with $k$ dominos appearing horizontally in row $1$.
\end{theorem}

\begin{proof}
    We first show $T\mapsto S_T$ is a bijection. Each element $T\in B_n$ is in natural bijection with its sequence of partial tableaux, $T_1,T_2,\dots T_n$, in which $T_i$ is $T$ with all but the dominos labeled $1,2,\dots,i$ removed: given a sequence take $T=T_n$ and given a tableaux, successively removing the highest labeled domino determines its partial tableaux sequence. We show that a given $S\subseteq [n]$ determines a unique partial tableaux sequence. Consider the insertions of the $i+1$ domino into $T_{i}$ for $i=1,2\dots,n-1$. 
    
    \textbf{Case 1: $i+1\in S$.} The $i+1$ domino must be placed horizontally in row $1$. 
    
    \textbf{Case 2: $i+1\not\in S$ and $\lambda_{i,1}\neq\lambda_{i,2}$.} The $i+1$ domino must be placed horizontally in row $2$, since a vertical placement would break left justification. Row lengths differ by an even number, so upper justification is respected.

    \textbf{Case 3: $i+1\not\in S$ and $\lambda_{i,1}=\lambda_{i,2}$.}
    The $i+1$ domino must be placed vertically, since a horizontal placement in row $2$ would not respect upper justification. 
    
   No matter the case, the $i+1$ domino has a unique placement in $T_{i}$ for all $i<n$. Consequentially, a given $S$ uniquely determines a sequence of partial tableaux and thus a unique $T$. This means the map $S\mapsto T$ such that $S_T=S$ is an injection. Injections between finite sets of equal cardinality are bijections, so $T$ is uniquely determined by the horizontal dominos appearing in row $1$. We find the desired binomial distribution by sending $T\mapsto |S_T|$, since there are $\binom{n}{k}$ sets $S_T\subseteq[n]$ with $k$ elements.
\end{proof}

\subsection{Rightward diagonal lattice walks}\label{subsec: lattice walks}

\begin{figure}[h]
    \centering
\begin{tikzpicture}
  \begin{axis}[xmin=0, xmax=2.5, ymin=-2, ymax=2, width=4cm, height=4cm, axis lines=middle, ticks=none]
    \addplot[mark=none] coordinates {(0,0) (1,-1) (2,-2)};
    \addplot[only marks, mark=*, mark size=2pt] coordinates {(1,-1) (2,-2)};
  \end{axis}
\end{tikzpicture}
\hspace{0.5cm}
\begin{tikzpicture}
  \begin{axis}[xmin=0, xmax=2.5, ymin=-2, ymax=2, width=4cm, height=4cm, axis lines=middle, ticks=none]
    \addplot[mark=none] coordinates {(0,0) (1,-1) (2,0)};
    \addplot[only marks, mark=*, mark size=2pt] coordinates {(1,-1) (2,0)};
  \end{axis}
\end{tikzpicture}
\hspace{0.5cm}
\begin{tikzpicture}
  \begin{axis}[xmin=0, xmax=2.5, ymin=-2, ymax=2, width=4cm, height=4cm, axis lines=middle, ticks=none]
    \addplot[mark=none] coordinates {(0,0) (1,1) (2,0)};
    \addplot[only marks, mark=*, mark size=2pt] coordinates {(1,1) (2,0)};
  \end{axis}
\end{tikzpicture}
\hspace{0.5cm}
\begin{tikzpicture}
  \begin{axis}[xmin=0, xmax=2.5, ymin=-2, ymax=2, width=4cm, height=4cm, axis lines=middle, ticks=none]
    \addplot[mark=none] coordinates {(0,0) (1,1) (2,2)};
    \addplot[only marks, mark=*, mark size=2pt] coordinates {(1,1) (2,2)};
  \end{axis}
\end{tikzpicture}

    \caption{The elements of $D_2$ with heights $0,0,1,2$, respectively}\label{fig: D2}
\end{figure}
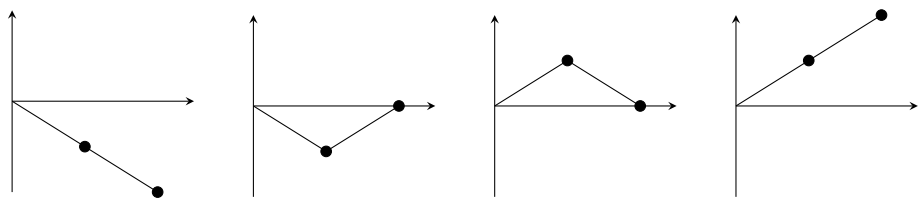

In \cite{ActonPetersenShirmanTenner2025}, the seemingly unrelated problem of counting $n$ step rightward diagonal lattice walks (RDLW) by their highest height emerged naturally in the study of a clairvoyant modification of Rob Pike's Malicious Ma\^ itre d' problem (see \cite{winkler2003}). Formally, a RDLW comprises successive steps from the origin to the upper right: $+(1,1)$ or lower right: $+(1,-1)$. Lattice walks are an important topic in Combinatorics, with many famous results, including the count of $n$ step RDLW terminating at point $\frac{\sqrt{2}}{2}(n+k,n-k)$ as $\binom{n+k}{k}$. We denote the set of $n$ step RDLW as $D_n$ and list $D_2$ in Figure $\ref{fig: D2}$.

At each step, we have a binary choice to move up or down, so $|D_n|=2^n$. With $|D_n|=|B_n|$, we may seek statistics on $D_n$ which mimic those of $B_n$. Representing elements $W\in D_n$ as binary strings of up steps (1) and down steps (0), we see the number of elements of $D_n$ with $k$ up steps is the number of ways to select $k$ entries from a size $n$ string to be $1$s. Then both are counted as $\binom{n}{k}$, giving the following result.

\begin{proposition}\label{prop: binomial RDLW stat}
    There are $\binom{n}{k}$ elements of $D_n$ with $k$ up steps.
\end{proposition}

Let the height of a RDLW be the maximum height of its vertical coordinates, and let $h:D_n\rightarrow\NN$ map $W\in D_n$ to its height. The authors of \cite{ActonPetersenShirmanTenner2025} show that this statistic gives the desired Pascalian distribution on elements of $D_n$.
\begin{theorem}\label{thm: pascalian RDLW stat}
    There are $\B{n}{k}$ elements of $D_n$ with height $k$.
\end{theorem}

\subsection{$B_n$ to $D_n$}\label{subsec: Bn to Cn}
The last result of this section describes a constructive bijection between $B_n$ and $D_n$ which unifies the binomial coefficient and Pascalian number statistics of each context. Given $T\in B_n$, construct $S_T$ and let step $n+1-s$ of $W\in D_n$ be up when $s\in S_T$ and down otherwise. The map $\phi$ which sends $T$ to this element of $D_n$ is the desired map, as we now show.
 
\begin{theorem}\label{Thm: Bn to Cn bijection}
    The map $\phi$ is a bijection from $B_n$ to $D_n$, under which
    \begin{enumerate}
        \item $T$ has $k$ horizontal dominos in its first row if and only if $\phi(T)$ has $k$ up steps, \label{first point}
        \item $T$ has columns of equal size if and only if  $\phi(T)$ has height $0$, and \label{second point}
        \item $T$ has shape $(n+k,n-k)$ if and only if $\phi(T)$ has height $k$.\label{third point}
    \end{enumerate}
    
\end{theorem}
\begin{proof}
Given $T$, let $S_T$ determine $D_n$'s up steps as described. Since elements of $D_n$ are uniquely determined by their up steps, $\phi$ is injective. Since $B_n$ and $D_n$ have equal size, $\phi$ is bijective. Thus, we need only show one direction of $\ref{first point},\ref{second point},$ and $\ref{third point}$. 

Because the size of $S_T$ is the number of up steps in $\phi(T)$, $\ref{first point}$ is immediate. Let $T_1,T_2,\dots T_n$ be $T$'s partial tableaux sequence (as described in Theorem $\ref{thm: ST determines T}$) and let $W_i=\phi(T_i)$. Let $T_i$ have shape $(\lambda_{i,1},\lambda_{i,2})$ and induce on $i$ for $\ref{second point}$ and $\ref{third point}$. To do this, we verify $\lambda_{i,2}+h(W_i)=n$.

When $i=1$, we have a single domino. If it lies vertically, $W_1$ is a single down step and has height $0$. If it lies horizontally, $W_1$ is a single up step with height $1$. In either case, $h(W_1)+\lambda_{1,2}=1$ and $h(W_1)=0$ if and only if $\lambda_{1,1}=\lambda_{1,2}$. Say inductively that $\lambda_{i,2}=n-h(W_i)$ and $\lambda_{i,1}=\lambda_{i,2}$ precisely when $h(\phi(T_i))=0$ for some $k<n$. We proceed by cases and note the placement of the $k+1$ domino is unique by Theorem $\ref{thm: ST determines T}$.

\textbf{Case $1$: $k+1\in S_T$.} Placing the $k+1$ domino horizontally in row $1$ corresponds to appending an up step at the start of $W_k$, so $h(W_{k+1})=1+h(W_{k})$. Since we've left row $2$ unchanged, $\lambda_{k+1,2}=\lambda_{k,2}=k-h(W_k)=(k+1)-h(W_{k+1})$ and $\ref{third point}$ holds. Further, $\ref{second point}$ holds, since $T_{k+1}$ has rows of unequal length and $h(W_{k+1})=1+h(W_{k})\geq 1$.

\textbf{Case $2$: $k+1\not\in S_T$ and $\lambda_{k,1}=\lambda_{k,2}$.} By the inductive hypothesis, $h(W_k)=0$. Placing the $k+1$ domino vertically corresponds to appending a down step at the start of $W_k$, so $h(W_{k+1})=h(W_{k})=0$. We've increased the size of both rows by $1$, so $\lambda_{k+1,2}=1+\lambda_{k,2}=1+k-h(W_k)=(k+1)-h(W_{k+1})$, giving $\ref{third point}$. Since $T_{k+1}$ has equal row lengths and $h(W_{k})=0$, $\ref{second point}$ holds.

\textbf{Case $3$ $k+1\not\in S_T$ and $\lambda_{k,1}>\lambda_{k,2}$.} By the inductive hypothesis, $h(W_k)>0$. Placing the $k+1$ domino horizontally in row $2$ corresponds to appending a down step at the start of $W_k$, so we've decreased its height by $1$. Then $h(W_{k+1})=h(W_{k})-1$. Since we increased row $2$'s length by $2$, we see $\lambda_{k+1,2}=2+\lambda_{k,2}=2+k-h(W_k)=(k+1)-h(W_{k+1})$, giving $\ref{third point}$. To see $\ref{second point}$ holds, say $\lambda_{k+1,1}=\lambda_{k+1,2}$, which forces each row to have length $k+1$. Since $k+1=\lambda_{k+1,2}+h(W_{k+1})$, we see $h(W_{k+1})$ must be zero. Alternatively, say $\lambda_{k+1,1}>\lambda_{k+1,2}$. Then $\lambda_{k+1,2}<k+1$ and $k+1=\lambda_{k+1,2}+h(W_{k+1})$ forces a positive $h(W_{k+1})$. Either way, $h(W_{k+1})=0$ precisely when $\lambda_{k+1,1}=\lambda_{k+1,2}$. 

Then no matter how the $k+1$ domino is placed, $\ref{first point},\ref{second point}$, and $\ref{third point}$ hold, and we are done.
\end{proof}

Theorem $\ref{Thm: Bn to Cn bijection}$ permits a natural interpretation of the Pascalian numbers on binary strings. Construct a size $n$ binary string $S$ for each $T\in B_n$ by letting entry $s$ of $S$ be $1$ when step $s$ of $\phi(T)$ is up and $0$ otherwise. Then define $h_s(S)$ to be the be number of $1$s at or before entry $s$ minus the number of $0$s at or before entry $s$. Then the largest amungst the $h_s$ is the height of $\phi(T)$.
\begin{corollary}
    There are $\B{n}{k}$ binary strings of length $n$ with $\max h_s=k$.
\end{corollary}

\section{Pascalian Polynomials}\label{sec: PP}
Any triangular array of numbers naturally defines to a sequence of polynomials. In our case we defije
$$R_n(x)=\sum_{k=0}^n\B{n}{k}x^k,$$
a sequence which starts $R_1(x)=1+x,R_2(x)=2+x+x^2,$ and $R_3(x)=3+3x+x^2+x^3$. In Subsection $\ref{subsec: polynomial recursions}$, we'll see the recursions for $\B{n}{k}$ derived in Section $\ref{sec: combinatorial foundations}$ extend to recursions on $R_n(x)$. We extract this recursion and formulate why $R_n(x)$ is best studied under the perspective of its reciprocal polynomial $P_n(x)=x^nR_n(1/x)$, which we call the $n$th Pascalian polynomial. We derive a nicer recursion for $P_n(x)$ and use it to derive generating functions for $P_n(x)$ and $R_n(x)$ in Subsection $\ref{subsec: GF for PP}$. 

\subsection{Recursions on the Pascalian Polynomials}\label{subsec: polynomial recursions} We begin with a recursion on $R_n(x)$.

\begin{proposition}\label{prop: Rn basic recursion}
    $R_{n}(x)=\left(x+\frac{1}{x}\right)R_{n-1}(x)+\B{n-1}{0}\left(1-\frac{1}{x}\right)$
\end{proposition}
\begin{proof}
We expand $R_{n}(x)$ as a sum in $\B{n}{k}$ and use the recursion from Theorem $\ref{thm: 2 col SDT satisfy pascalian recursion}$ for terms of degree $0<k<n-1$ to see
\begin{align*}
    R_{n}(z)=&\sum_{k=0}^{n}\B{n}{k}x^k=\B{n}{0}+x^{n-1}+x^{n}+\sum_{k=1}^{n-2}\B{n}{k}x^k\\
    =&\B{n}{0}+x^{n-1}+x^{n}+\left(\sum_{k=1}^{n-2}\B{n-1}{k+1}x^k\right)+\left(\sum_{k=1}^{n-2}\B{n-1}{k-1}x^k\right)\\
    =&\B{n}{0}+x^{n-1}+x^{n}+\frac{1}{x}\left(R_{n-1}(x)-\B{n-1}{0}-\B{n-1}{1}x\right)\\
    +&x\left(R_{n-1}(x)-x^{n-2}-x^{n-1}\right).
    \end{align*}
Collecting terms and recalling $\B{n+1}{0}-\B{n}{1}=\B{n}{0}$ gives the desired result.
\end{proof}
We now consider the recursion for $P_n(x)$, which turns out to be much easier to work with. By mapping $x\mapsto 1/x$ and multiplying through by $x^n$ in $R_n(x)$'s recursion from Proposition $\ref{prop: Rn basic recursion}$ yields the following result.

\begin{proposition}\label{prop: Pn Recursion for k=0}
   $P_{n}(x)=(1+x^2)P_{n-1}(x)+\B{n-1}{0}(1-x)x^{n}$
\end{proposition}
We may use Proposition $\ref{prop: Pn Recursion for k=0}$ to see that no consecutive pair of Pascalian polynomials may share a root: if $P_n(a)$ and $P_{n-1}(a)$ vanish, then $\B{n-1}{0}(1-a)a^{n}=0$ and $a=0$ or $a=1$--impossible since neither are roots of $P_n(x)$ for any $n$. Successive application of this recursion, done easily with no rational functions in sight, will express $P_n(x)$ in terms of $P_{n-k}(x)$. A simple induction gives the desired expression.

\begin{theorem}\label{thm: extended Pn recursion, general k}
    For all $k=1\dots n-1$ we have $$P_{n}(x)=(1+x^2)^{k}P_{n-k}(x)+(1-x)\sum_{j=0}^{k-1}\B{n-j-1}{0}(1+x^2)^jx^{n-j}.$$
\end{theorem}

\begin{proof}
    If $k=1$, we recover Proposition \ref{prop: Pn Recursion for k=0}, so assume our claim for some $k<n$. Substituting $(1+x^2)P_{n-k-1}(x)+\B{n-k-1}{0}(1-x)x^{n-k}$ for $P_{n-k}(x)$ into the inductive hypothesis and distributing shows
    \begin{align}
        P_{n}(x)&=(1+x^2)^{k}P_{n-k}(x)+(1-x)\sum_{j=0}^{k-1}\B{n-j-1}{0}(1+x^2)^jx^{n-j}\nonumber\\
        &=(1+x^2)^{k+1}P_{n-k-1}(x)+ \nonumber\\
        &\B{n-k-1}{0}(1+x^2)^{k}(1-x)x^{n-k}+(1-x)\sum_{j=0}^{k-1}\B{n-j}{0}(1+x^2)^jx^{n+1-j}.\label{eq: repeated recursion}
    \end{align}
    The terms in $\eqref{eq: repeated recursion}$ combine into the single sum $(1-x)\sum_{j=0}^{k}\B{n-j}{0}(1+x^2)^jx^{n+1-j},$
    completing the desired inductive step. 
\end{proof}
Comparing coefficients above yields a ``higher order'' recursion on $\B{n}{k}$, equivalent to successive application of Theorem $\ref{thm: 2 col SDT satisfy pascalian recursion}$. The statement written equivalently with the numbers $\B{n}{k}$ themselves splits into many unruly cases based on $k$'s distance to $n$. Combinatorially, Theorem $\ref{thm: extended Pn recursion, general k}$ nicely packages the precise construction of elements in $B(n,k)$ by concatenating the last $m\leq n$ dominos onto smaller tableaux.

\subsection{Generating Functionology}\label{subsec: GF for PP}
We now derive the generating functions for $P_n(x)$ and $R_n(x)$ using a well-known generating function for the numbers $\B{n}{0}$ from \cite[A001405]{OEIS}. Particularly, for $|z|<1/2$,
$$\sum_{n\geq0}\B{n}{0}z^n=\frac{1}{2z}\left(\sqrt{\frac{1+2z}{1-2z}}-1\right)$$ 
\begin{theorem}\label{thm: GF for Pn}
    The power series $\sum_{n\geq0}P_n(x)z^n$ and $\sum_{n\geq0}R_n(x)z^n$ have respective generating functions 
    
    $$G(x,z)=\frac{2+(x-1)\left(1-\sqrt{\frac{1+2xz}{1-2xz}}\right)}{2\left(1-z(1+x^2)\right)}\hspace{.5in}H(x,z)=\frac{2x+(1-x)\left(1-\sqrt{\frac{1+2z}{1-2z}}\right)}{2(x-z(1+x^2))}.$$ $G(x,z)$ converges when $|z(1+x^2)|<1$ and $|xz|<1/2$ while $H(x,z)$ converges when $|z(1+x^2)|<|x|$ and $|2z|<1$. 
\end{theorem}
\begin{proof}
We first expand $G(x,z)$ and use Proposition $\ref{prop: Pn Recursion for k=0}$ to write $P_n(x)$ in terms of $P_{n-1}(x)$:
    \begin{align*}
        G(x,z)=&1+\sum_{n\geq1}P_n(x)z^n\\
        =&1+z\sum_{n\geq0}\left((1+x^2)P_n(x)+\B{n}{0}(1-x)x^{n+1}\right)z^n\\
        =&1+z(1+x^2)G(x,z)+xz(1-x)\sum_{n\geq0}\B{n}{0}(xz)^n\\
        =&1+z(1+x^2)G(x,z)+\frac{1}{2}(x-1)\left(1-\sqrt{\frac{1+2xz}{1-2xz}}\right)\\
    \end{align*}
    Solving this as a linear equation in $G(x,z)$ yields the desired expression. Because $P_n(z)$ and $R_n(z)$ are reciprocal polynomials, we recover $H(x,z)$ as $G\left(\frac{1}{x},xz\right)$. Making the desired substitution in $G(x,z)$ and clearing denominators gives the desired result
\end{proof}
The explicit relation $P_n(x)+x^{n+1}R_n(x)=(1+x)(1+x^2)^n$, which will be absolutely vital in Section \ref{sec: limit curve of roots}, can be viewed as a relationship between $G(x,z)$ and $H(x,z)$. The coefficient of $z^n$ in $xH(x,xz)$ is $x^{n+1}R_n(x)$, so $$G(x,z)+xH(x,xz)=\frac{1+x}{1-z(1+x^2)}=\sum_{n\geq0}(1+x)(1+x^2)^nz^n.$$

We conclude this section with an additional recursion for $P_n(x)$ that follows from isolating the radical in $G(x,z)$, squaring both sides, and collecting terms. 
\begin{theorem}\label{thm: PP linear decomposition}
    $P_n(x)$ satisfies the recursion $$P_n(x)=(2x)^n+(1-x)\sum_{j=0}^{n-1} \B{n-j-1}{0}P_j(x)x^{n-j-1}.$$
\end{theorem}
\begin{proof}
Isolating the root in the generating function and squaring gives
\begin{align*}
    \sqrt{\frac{1+2xz}{1-2xz}}=&\frac{1}{(1-x)}(2G(x,z)(1-z(1+x^2))-(1+x))\implies\\
    \frac{1+2xz}{1-2xz}=&\frac{1}{(1-x)^2}\Big(4G(x,z)^2(1-z(1+x^2))^2-4(1+x)G(x,z)(1-z(1+x^2))+(1+x)^2\Big).\\
\end{align*}
Moving terms to one side and multiplying through by $(1-x)^2$ yields
$$0=4G(x,z)^2(1-z(1+x^2))^2-4G(x,z)(1-z(1+x^2))(1+x)+(1+x)^2-\frac{(1-x)^2(1+2xz)}{1-2xz}.$$
However, the constant term in $G(x,z)$ combines as
$$(1+x)^2-\frac{(1-x)^2(1+2xz)}{1-2xz}=\frac{4x\big(1-z(1+x^2)\big)}{1-2xz}.$$
Substituting this and dividing a common factor of $4(1-z(1+x^2))$ from the resulting equation gives the form
\begin{equation}
    0=G(x,z)^2(1-z(1+x^2))-G(x,z)(1+x)+\frac{x}{1-2xz}\makelabel\label{eq: vanishing coeff of z^n}
\end{equation}
The right side of equation $\eqref{eq: vanishing coeff of z^n}$ may be written as a power series in $z$, which must have all vanishing coefficients. We now isolate the coefficient of $z^n$ in each term above, starting with the constant term in $G(x,z):$
    $$\frac{x}{1-2xz}=x\sum_{n\geq0}(2x)^nz^n=\sum_{n\geq0}2^{n}x^{n+1}z^n.$$
Thus we find coefficients in $z$ of $2^{n}x^{n+1}$. The coefficient of $z^n$ in $(1+x)G(x,z)$ is $(1+x)P_n(x)$, so we examine $(1-z(1+x^2))G(x,z)^2$. This expands as
$$\sum_{n\geq0}\left(\sum_{j=0}^nP_j(x)P_{n-j}(x)\right)z^n-(1+x^2)\sum_{n\geq0}\left(\sum_{j=0}^nP_j(x)P_{n-j}(x)\right)z^{n+1}.$$
Peeling the $n=0$ term from the first $n$ sum and reindexing the second gives the form
$$1+\sum_{n\geq1}\left(\sum_{j=0}^nP_j(x)P_{n-j}(x)\right)z^n-(1+x^2)\sum_{n\geq1}\left(\sum_{j=0}^{n-1}P_j(x)P_{n-1-j}(x)\right)z^{n}.$$
Now we pull the $j=n$ term from the first inner sum and combine sums to find
$$1+\sum_{n\geq1}\left(P_n(x)+\sum_{j=0}^{n-1} P_j(x)\left(P_{n-j}(x)-(1+x^2)P_{n-1-j}(x)\right)\right)z^{n}.
$$
Proposition $\ref{prop: Pn Recursion for k=0}$ tells us $P_{n-j}(x)-(1+x^2)P_{n-1-j}(x)=\B{n-j-1}{0}x^{n-j}(1-x)$ for all $j=0,1,\dots n-1$, so our expression for $(1-z(1+x^2))G(x,z)^2$ reduces to
    $$1+\sum_{n\geq1}\left(P_n(x)+\sum_{j=0}^{n-1} \B{n-j-1}{0}P_j(x)(1-x)x^{n-j}\right)z^{n},$$
with coefficients in $z$ equal to $P_n(x)+\sum_{j=0}^{n-1} \B{n-j-1}{0}P_j(x)(1-x)x^{n-j}$. The final coefficient of $z^n$ in $\eqref{eq: vanishing coeff of z^n}$ vanishes for all $n$, so we see
$$0=P_n(x)+\sum_{j=0}^{n-1}\left(\B{n-j-1}{0}P_j(x)(1-x)x^{n-j}\right)-(1+x)P_n(x)+x(2x)^n.$$
The desired result then follows from canceling $\pm P_n(x)$ and dividing $x$ from the remaining terms. 
\end{proof}

\section{Roots of Pascalian Polynomials}\label{sec: roots of PP}

Here we pass to a complex variable $z$ and study the roots of Pascalian Polynomials. We start with bounds in the complex plane before counting the rational, real, and purely imaginary roots of $P_n(z)$. It will be convenient to introduce the polynomial
    $$q_n(z)=\sum_{k=0}^{\lfloor n/2\rfloor}\binom{n}{k}z^k, \hspace{.4in} \text{so}\hspace{.4in} P_n(z)=(1+z)q_n(z^2)+\begin{cases}
        \B{n}{0}z^n, \,\,\,\,&n \text{ even}\\
        0, \,\,\,\,&n \text{ odd}.\\
    \end{cases}$$
The polynomial $q_n(z)$ belongs to a larger family of truncated binomial expansions, which has seen expansive interest in the past two decades. Active work on these polynomials will underlie key results concerning Pascalian polynomials. Scherbak initiated the study in the $2004$ MSRI program on topological aspects of real algebraic geometry. She showed the truncated binomial polynomials formed a natural basis for key intersections of certain Schubert varieties, thus giving deep motivation for a study for their irreducibility, examined in \cite{DubickasSiurys2017,FilasetaKumchevPasechnik2007,KhandujaKhassaLaishram2011,laishram-yadav2024}. 
\subsection{Bounds on the roots of $P_n(z)$ within the complex plane}\label{subsec: location of PP roots in complex plane}
\begin{figure}[htbp] 
  \centering
  \begin{minipage}[t]{0.3\linewidth}
    \centering
    \includegraphics[width=\linewidth]{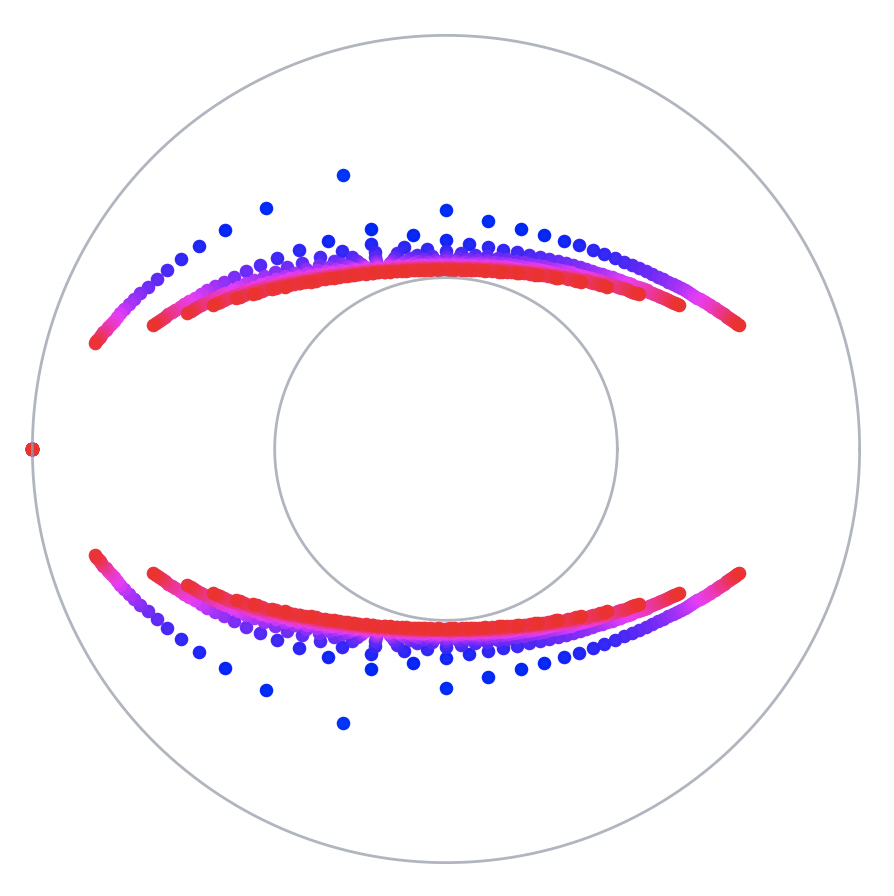}
  \end{minipage}
  \caption{The roots of $P_n(z)$ plotted for $n\leq50$ with the annulus of Theorem $\ref{thm: PP roots annuli}$ in gray. The color of $P_n(z)$'s roots flows from blue to red as $n$ grows.}
  \label{fig: roots in annulus} 
\end{figure}

By writing $P_n(z)$ in terms of $q_n(z)$, we see $P_n(-1)$ is zero for odd $n$ and $\B{n}{0}$ for even $n$. It will be convenient to call this root at $-1$ the \textit{trivial} Pascalian root. We now state and prove the optimal annulus which bounds the roots of $P_n(z)$, as visualized in Figure $\ref{fig: roots in annulus}$. 

\begin{theorem}\label{thm: PP roots annuli}
    The nontrivial roots of $P_n(z)$ lie within the annulus $\sqrt{2}-1<|z|<1$.
\end{theorem}

\begin{proof}
     We begin with the lower bound, so set $r=\sqrt{2}-1$ and consider $|z|<r$. As noted after Theorem $\ref{thm: GF for Pn}$, $(1+z)(1+z^2)^n=P_n(z)+z^{n+1}R_n(z)$. Solving for $P_n(z)$ and dividing by $(1+z)(1+z^2)^n$ shows
    $$\left|\frac{P_n(z)}{(1+z)(1+z^2)^n}\right|\geq1-\left|\frac{z^{n+1}R_n(z)}{(1+z)(1+z^2)^n}\right|= 1-\left|\frac{z}{1+z}\right|\left|\frac{z}{1+z^2}\right|^n|R_n(z)|.$$
    We bound $|P_n(z)|$ strictly above $0$ by bounding the final term on the right above by $1$. By the maximum modulus principal, its maximum occurs along $|z|=r$. On this contour, the numerators of $\frac{z}{1+z}$ and $\frac{z}{1+z^2}$ are both fixed, so each factor is maximized where its denominator is minimized. The first then achieves its maximum of $\sqrt{2}/2$ on the negative real axis, while the second achieves its maximum of $1/2$ on the imaginary axis. As these maxima are achieved at distinct $z$, we have the strict inequality
    $$\left|\frac{z}{1+z}\right|\left|\frac{z}{1+z^2}\right|^n|R_n(z)|<\frac{\sqrt{2}}{2^{n+1}}|R_n(z)|.$$
    As a polynomial with positive real coefficients, $|R_n(z)|$ is maximized on the positive real axis. We factor out the constant term and find
    $$|R(z)|=\B{n}{0}\left|\sum_{k=0}^n\B{n}{k}\B{n}{0}^{-1}z^{k}\right|< \B{n}{0}\left(\frac{1-z^{n+1}}{1-z}\right)\leq\B{n}{0}\left(\frac{1-r^{n+1}}{1-r}\right),$$
    where the second step follows from bounding the polynomial's coefficients above by $1$. Since $\sqrt{2}<3/2$ and $(1-r)^{-1}<2$, we observe
    $$\frac{\sqrt2}{2^{n+1}}\left(\frac{1-r^{n+1}}{1-r}\right)\leq \frac{3}{2^{n+1}}.$$ 
    Finally, we see $\sum_{k=0}^{\lfloor n/2\rfloor}\binom{n}{k}\geq 2^{n-1}$, so $\binom{n}{\lfloor n/2 \rfloor}\leq 2^{n-1}$. Combining these results shows
    $$\left|\frac{z}{1+z}\right|\left|\frac{z}{1+z^2}\right|^n|R_n(z)|<3/4.$$
    As for the upper bound, $P_n(z)$ has weakly increasing coefficients, so the Eneström-Kakeya theorem (\cite[Thm. 1.1.5]{Prasolov2001} bounds the roots within the closed unit disk. For $n$ odd, $P_n(z)$'s nontrivial roots as those of $q_n(z^2)$. Since $q_n(z)$ has strictly increasing coefficients, its roots, and those of $q_n(z^2)$, are bound within the open unit disk by the Eneström-Kakeya theorem. We now verify the remaining roots lie strictly within the disk $|z|\leq1$, so consider $n=2m$. With $P_{2m}(1)=2^{n-1}$, we multiply

    $$(1-z)P_{2m}(z)=1+\sum_{k=1}^{2m}\left(\B{n}{n-k}-\B{n}{n-k+1}\right)z^k-\B{2m}{0}z^{2m+1}$$
    and note $\B{n}{n-k}-\B{n}{n-k+1}$ is zero for odd $k$ and positive for even $k$. This lets us write 

    $$U_{2m}(z)=1+\sum_{k=1}^m\left(\B{2m}{2m-2k}-\B{2m}{2m-2k+1}\right)z^k,$$
    $$\text{so} \hspace{.2in}(1-z)P_{2m}(z)=U_{2m}(z^2)-\B{2m}{0}z^{2m+1}.$$
    We note $U_{2m}(z)$ has strictly positive coefficients with $U_{2m}(1)$ telescoping to $\B{2m}{0}$. Thus if $P_{2m}(z)=0$ and $|z|=1$, we must have $|U_{2m}(z^2)|=\B{2m}{0}$. However, the triangle inequality bounds $|U_{2m}(z^2)|$ for $z\neq\pm1$ strictly below its value at $z=\pm1$, which is precisely $\B{2m}{0}$. Then the solutions of $P_{2m}(z)=0$ along $|z|=1$ must be $z=\pm1$, neither of which is a root of $P_{2m}(z)=0$. With this, we are done.
\end{proof}
Because the roots of $R_n(z)$ are the reciprocals of $P_n(z)$'s roots, we have the immediate corollary for $R_n(z)$'s roots.

\begin{corollary}\label{thm: Rn roots annulus}
    The nontrivial roots of $R_n(z)$ lie within the annulus $1< |z|<1+\sqrt{2}$.
\end{corollary}

Since all but at most one of $P_n(z)$'s roots has norm less than $1$, their product goes to $0$ as $n$ grows. Vietta's relations affirm this with the stronger result that the product of $P_n(z)$'s roots is precisely $1/\binom{n}{\lfloor n/2\rfloor}$.

In the next section, we'll show roots of $P_n(z)$ may be found arbitrarily close to the points $\pm i(\sqrt{2}-1)$ and $\pm1$, so both bounds are optimal.  This is illustrated in Figure $\ref{fig: roots in annulus}$.

Theorem $\ref{thm: PP roots annuli}$ yields further observations for the roots of $P_n(z)$. Particularly, the strict bound of $P_n(z)$'s nontrivial roots within the open unit disk lets us classify the common roots of $P_n(z)$ and $P_{n-2}(z)$.

\begin{corollary}\label{prop: F(n,2)(z)}
    $P_n(z)$ and $P_{n-2}(z)$ share no nontrivial roots for all $n\geq2$.
\end{corollary}
\begin{proof}
   Setting $k=2$ in Theorem $\ref{thm: extended Pn recursion, general k}$ shows 
    $$P_n(z)=(1+z^2)^2P_{n-2}(z)+(1-z)z^{n-1}\left(\B{n-1}{0}z+\B{n-2}{0}(1+z^2)\right)$$
Particularly, if $P_n(z)=P_{n-2}(z)=0$ for some $z$, then the remaining term
\begin{equation}
    0=(1-z)z^{n-1}\left(\B{n-1}{0}z+\B{n-2}{0}(1+z^2)\right)\label{eq: extra term}
\end{equation}
also vanishes. If $n$ is odd, then $\B{n-1}{0}=2\B{n-2}{0}$, and the expression above factors as $(1-z)z^{n-1}\B{n-2}{0}(1+z)^2$, with roots at $-1, 0,$ and $1$. Since no Pascalian polynomial has $1$ or $0$ as a root, $P_n(z)$ and $P_{n-2}(z)$ may only share the trivial root for $n$ odd (and in fact, do).

For $n=2s$, we see $\B{n-1}{0}z+\B{n-2}{0}(1+z^2)$ factors as $$\cat_{s-1}(s+(2s-1)z+sz^2),$$ where $\cat_n=\binom{2n}{n}/(n+1)$ denotes the $n$th Catalan number \cite[A000108]{OEIS}. Expression $\eqref{eq: extra term}$ then has roots $1,0$, or $\frac{1}{2s}\left(2s-1\pm i\sqrt{4s-1}\right)$. No pascalian polynomial vanishes at $0$, and the remaining options for $z$ have norm exactly equal to $1$. Since all roots of $P_n(z)$ have norm strictly less than $1$ for even $n$, no such $z$ simultaneously satisfies $P_n(z)=P_{n-2}(z)=0$ and expression $\eqref{eq: extra term}$. Thus $P_n(z)$ and $P_{n-2}(z)$ cannot share a common root.
\end{proof}

We generalize this in the following conjecture.

\begin{conjecture}\label{conj: common roots of Pn(z) and Pn-k(z)}
    $P_n(z)$ and $P_{n-k}(z)$ share no nontrivial roots for all $n$ and $1\leq k\leq n$.
\end{conjecture}

\subsection{Real and imaginary roots}\label{subsec: real, imagionary roots of PP}
We now simultaneously classify the integral, rational, and real roots of $P_n(z)$ as precisely the trivial at $-1$ for odd $n$.
\begin{lemma}\label{lem: unique real root}
    The trivial root at $-1$ is the unique real root of $P_n(z)$ for odd $n$. For even $n$ and real $z$, $P_n(z)$ is strictly positive.
\end{lemma}
\begin{proof}    
    Consider odd $n$, so $P_n(z)=(1+z)q_n(z^2)$. Because $q_n(z^2)$ is an even function with positive coefficients, its global minima is its constant term $q_n(0)=1>0$. Since $q_n(z^2)$ is strictly positive and $q_n(z^2)$'s roots are $P_n(z)$'s nontrivial roots, $P_n(z)$ has no other real roots.

    Now consider even $n$, so $P_n(z)=(1+z)q_n(z^2)+\binom{n}{n/2}z^n$. $P_n(z)$ has all positive coefficients and nonzero constant term, so its real roots are negative. Using Theorem $\ref{thm: PP roots annuli}$, we limit our view to the interval $(-1,0)$. On this interval both $(1+z)q_n(z^2)$ and $\binom{n}{n/2}z^n$ are strictly positive, so we are done.
\end{proof}
This result sheds light on sum of $P_n(z)$'s roots roots in $\CC$, which we view as their ``center of mass.'' Particularly, if $z_0$ is a nontrivial root of $P_n(z)$, then so are its conjugate, negative, and negative conjugate. These individual sets have sum $0$, so the sum of $P_n(z)$'s roots is $-1$ for odd $n$. To understand the roots' sum for even $n$, we appeal to Vietta's relation and find they add to $-\B{n}{1}\B{n}{0}^{-1}=2/(n+2)-1$, which approaches $-1$. Turning to theseth roots in $\CC$, our examination of the purely imaginary case is more complicated.

\begin{theorem}\label{prop: no imag root for even n}
    $P_n(z)$ has no purely imaginary roots for $n\equiv 0,1,2\mod 4$. For $n\equiv3\mod4$, $P_n(z)$ has a unique pair of purely imaginary roots.
\end{theorem}
\begin{proof}
\textbf{Case 1: $n=2m$ even.} We separate $P_n(xi)$ into its real and imaginary parts for arbitrary real $x$ as
$$\sum_{k=0}^{n}\B{n}{n-k}(xi)^k=\sum_{k=0}^{m}(-1)^k\binom{n}{k}x^{2k}+ix\sum_{k=0}^{m-1}(-1)^k\binom{n}{k}x^{2k}.$$
When the real part of the above expression vanishes, the imaginary part is precisely $(-1)^{m}\binom{n}{m-1}x^{n}$. This must also vanish, so $x=0$. Since $P_n(0)=1$, no such root $xi$ may exist. 

\textbf{Case 2: $n\equiv 1\mod4$.} Since the imaginary roots of $P_n(x)$ the are precisely $i$ times the real roots of $q_n(-x^2)$, we prove $q_n(-x^2)$ has no roots on $(-1,1)$. It suffices to show $q_n(-x)$ has no roots on $(0,1)$. To do so, set $n=2m+1$ and recognize
\begin{equation}
q_n(-x)=(1-x)^n+\sum_{k=m+1}^n(-1)^{k+1}\binom{n}{k}x^k.\makelabel\label{eq: qn(-x) tail}
\end{equation}
Peeling the $k=n$ term from the second term of Equation $\eqref{eq: qn(-x) tail}$ and reindexing with $k\mapsto m+k$ shows 
\begin{align*}
    &x^n+\sum_{k=1}^{m}(-1)^{k-1}\binom{n}{m+k}x^{m+k}\\
    =&x^n+\sum_{k=1}^{\frac{m}{2}}\binom{n}{m+2k-1}x^{m+2k-1}-\binom{n}{m+2k}x^{m+2k}.
\end{align*}
Each term of the above sum is $x^{m+2k-1}(\binom{n}{m+2k-1}-\binom{n}{m+2k}x)$, which is strictly positive since $\binom{n}{m+2k-1}>\binom{n}{m+2k}$ and $x\in(0,1)$. Equation $\eqref{eq: qn(-x) tail}$ then expresses $q_n(-x)$ as a sum of strictly positive expressions and consequentially guarantees $q_n(-x)$'s positivity on $(0,1)$. 

\textbf{Case 3: $n\equiv3\mod4$}. As before, set $n=2m+1$. We differentiate $q_n(-x)$, reindex with $k\mapsto k-1$ substitute $n\binom{n-1}{k-1}=k\binom{n}{k}$ to find
$$\frac{d}{dx}q_n(-x)=\sum_{k=1}^{m}(-1)^k\binom{n}{k}kx^{k-1}=-n\sum_{k=0}^{m-1}(-1)^k\binom{n-1}{k}x^{k}.$$
Similarly to the proceeding case, we expand $-\frac{1}{n}\frac{d}{dx}q_n(-x)$ as
\begin{equation}
    \sum_{k=0}^{m-1}(-1)^k\binom{n-1}{k}x^{k}=(1-x)^{n-1}+\sum_{k=m}^{n-1}(-1)^{k-1}\binom{n-1}{k}x^{k}\makelabel\label{eq: qn imag roots 3mod4}.
\end{equation}
Reindexing $k\mapsto k+m$ in the rightmost sum above and pairing successive terms yields
\begin{align*}
    &\sum_{k=0}^{m+1}(-1)^{k}\binom{n-1}{m+k}x^{m+k}\\
    =&\sum_{k=0}^{\frac{m+1}{2}}\binom{n-1}{m+2k}x^{m+2k}-\binom{n-1}{m+2k+1}x^{m+2k+1}.
\end{align*}
Each term of the above sum is $x^{m+2k}(\binom{n-1}{m+2k}-\binom{n-1}{m+2k+1}x)$, which is strictly positive since $\binom{n-1}{m+2k-1}>\binom{n-1}{m+2k+1}$ and $x\in(0,1)$. Equation $\eqref{eq: qn imag roots 3mod4}$ then expresses $-\frac{1}{n}\frac{d}{dx}q_n(-x)$ as a sum of strictly positive expressions, so $q_n(-x)$ has a strictly negative derivative. Then evaluating $q_n(-1)$ with the truncated alternating binomial coefficient sum 
    $$\sum_{k=0}^D(-1)^k\binom{n}{k}=(-1)^{D}\binom{n-1}{D},\;\;\;\,\,\,\text{ we see }\;\;\;\,\,\, q_n(-1)=(-1)^{m}\binom{n-1}{m}.$$
Since $n\equiv3\mod4$, this shows $q_n(-1)<0$. Because $q_n(0)=1$, $q_n(-1)=-\binom{n}{m}$, and $q_n(-x)$ is strictly decreasing, the intermediate value theorem guarantees a unique real root of $q_n(-x)$ on $(0,1)$, which extends to a real root pair of $q_n(-x^2)$. Thus, $P_n(z)$ has a unique conjugate pair of purely imaginary roots if and only if $n\equiv 3\mod 4$.

\end{proof}

\section{Convergence of Roots to $\Gamma_n$}\label{sec: limit curve of roots}
The asymptotic behavior of the roots of  a family of polynomials is a well-known problem of Complex Function Theory. Jentzsch \cite{Jentzsch1917} initiated the study for families which are a truncations of power series with finite radii of convergence, and Szeg\H{o} \cite{szego1924} famously investigated the root asymptotics of normalized truncations of the Taylor series of $e^z$--giving the Szeg\H{o} curve $|ze^{1-z}|=1$ which the roots tend to asymptotically and fill densely. Since then, the general problem has expanded greatly with an explosion of methods. 

While popular and classical methods rely heavily on approximating integral representations with Laplace's method and Watson's Lemma (see \cite{deBruijn1981}), we employ more elementary means. Our arguments are largely inspired by those of Janson's and Norfolk's 2009 paper \cite{janson2009} which  studied the root asymptotics of the truncated binomial polynomials. In this section, we give a family of curves $\Gamma_n$ which bounds and well approximates the roots of $P_n(z)$ and which asymptotically accumulates the roots densely. 

\subsection{Properties of Key Curves}\label{subsec: properties of key curves} 
\begin{figure}[htbp]
  \centering
  \hspace*{.7in}
  \begin{minipage}[t]{0.22\linewidth}
    \centering
    \includegraphics[width=\linewidth]{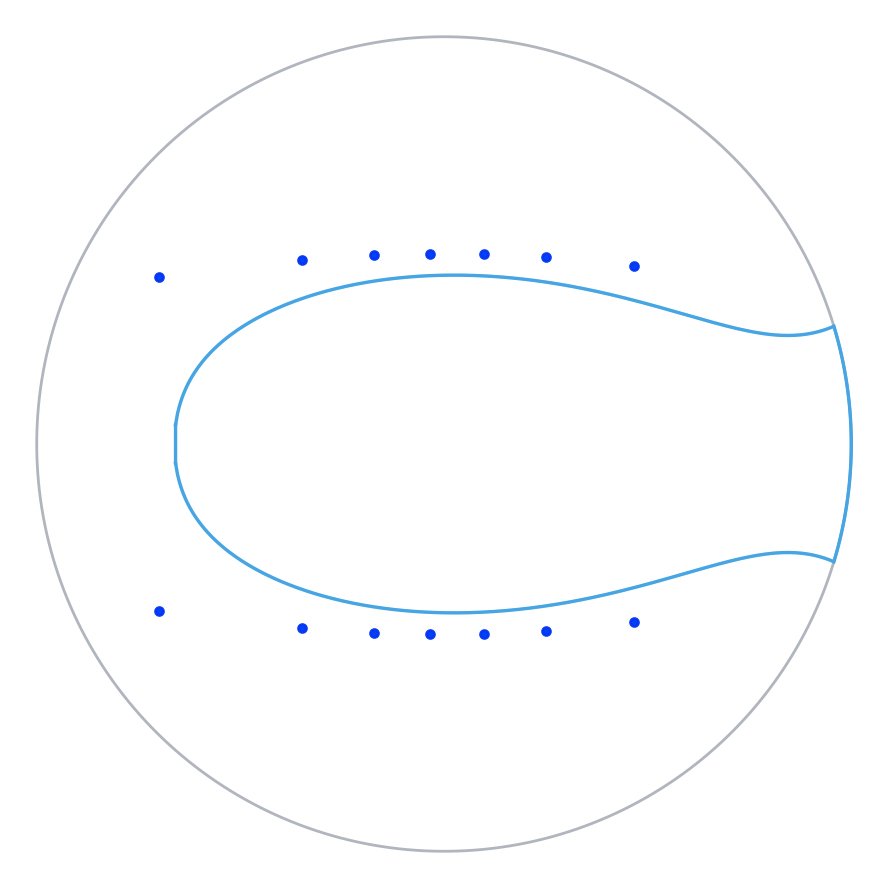}
  \end{minipage}\hfill
  \begin{minipage}[t]{0.22\linewidth}
    \centering
    \includegraphics[width=\linewidth]{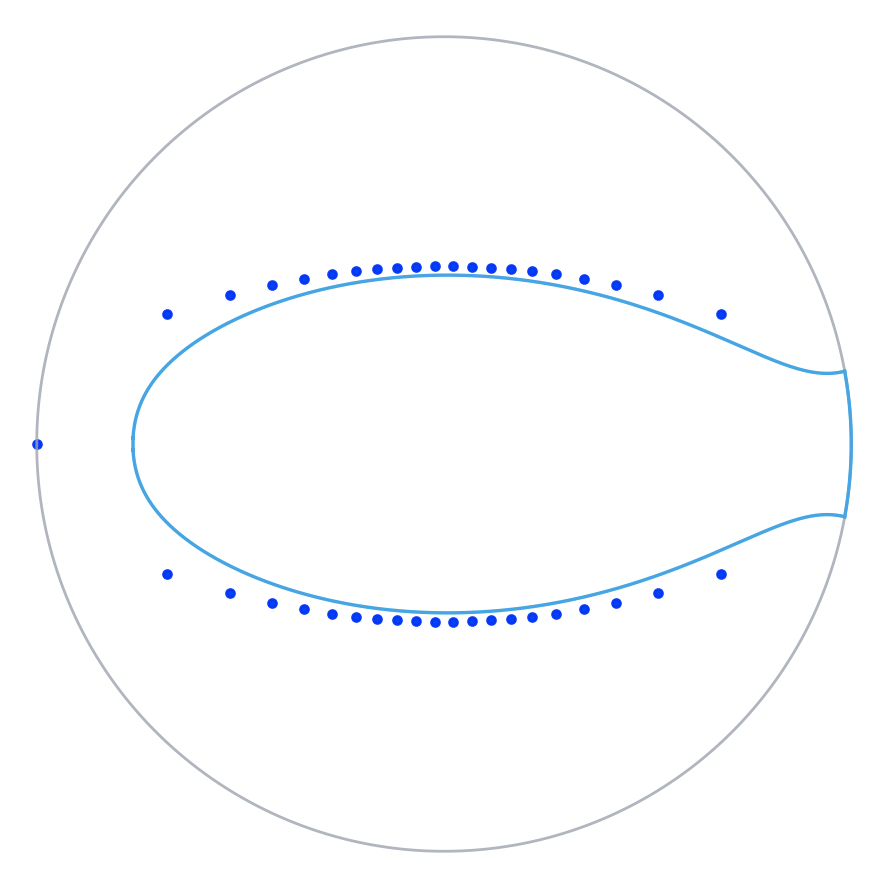}
  \end{minipage}\hfill
  \begin{minipage}[t]{0.22\linewidth}
    \centering
    \includegraphics[width=\linewidth]{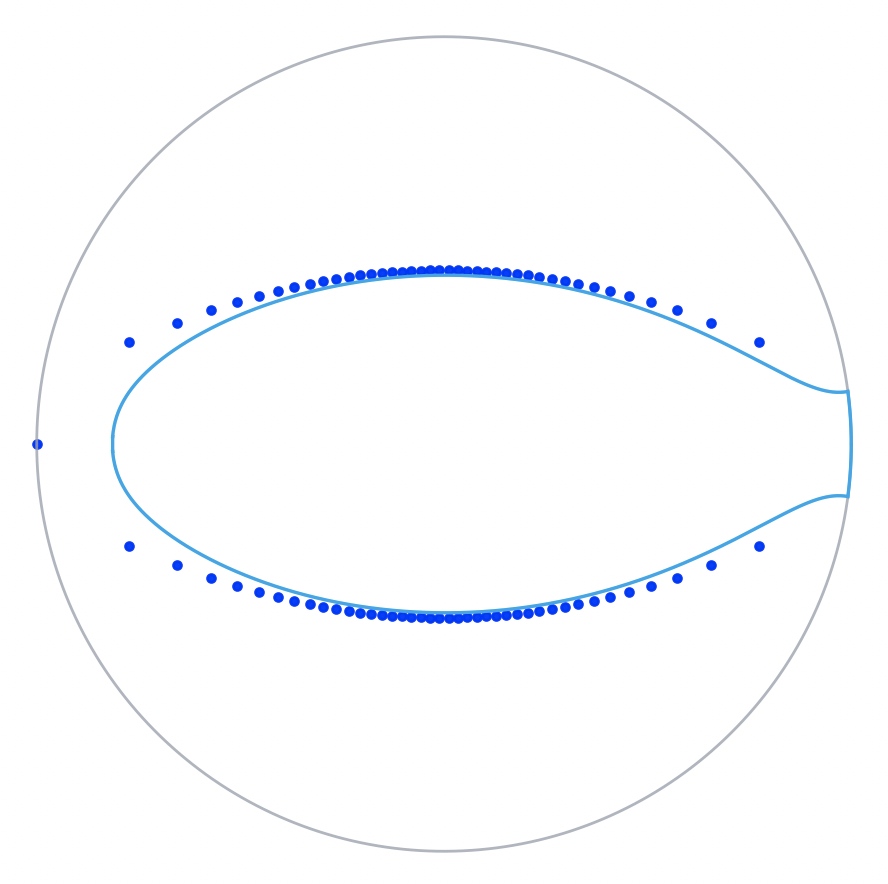}
  \end{minipage}\hspace*{.7in}
  \caption{The roots of $P_n(z)$ bound by $\partial\Gamma_n$ for $n=14,41,81$.}
  \label{fig: Pn roots bound by Cn}
\end{figure}

We start by defining the family of curves we use. Set $n\geq2$,
$$\frac{3}{8}\leq K_n=\frac{n^2-1}{2n^2}\leq\frac{1}{2},$$
$$\Gamma_n = \left\{ z\in\CC:\ \frac{|z|}{\sqrt[n]{|1+z|}|1+z^2|}\leq K_n,\ |z|\le 1\right\},$$
and $\partial\Gamma_n$ to the boundary of $\Gamma_n$. See $\ref{fig: Pn roots bound by Cn}$ for examples of these curves for select $n$. We now prove the following lemma about $\Gamma_n$'s asymptotics. We let $\Gamma$ denote the limit of $\Gamma_n$ and $\partial\Gamma$ is interpreted similarly.
\begin{lemma}\label{lem: propoerties of Cn}
    For all finite $n$, $\partial\Gamma_n$ contains the point $z=1$ and $\Gamma_n$
    converges to the intersection of the disks centered at $\pm i$ with radii $\sqrt{2}$. 
\end{lemma}

\begin{proof}
    First, we verify that $z=1$ lies in $\Gamma_n$. At $z=1$, $$\frac{|z|}{\sqrt[n]{|1+z|}|1+z^2|}=\frac{1}{2^{1+1/n}}\leq \frac{\sqrt{2}}{4}<\frac{3}{8}\leq K_n$$
    where the first inequality minimizes $2^{-(n+1)/n}$ at $n=2$ and the second is bounding $\sqrt{2}<1.5$. Since $|z|\leq 1$ on $\Gamma_n$, the point $z=1\in \Gamma_n$ must belong to the boundary $\partial\Gamma_n$.

    Now we describe $\Gamma$, so consider $z\in\Gamma_n$ and thus satisfies $\frac{|z|}{|z^2+1|}\leq \frac{1}{2}$. Clearing denominations, squaring, and rearranging shows $0\leq|z^2+1|^2-4|z|^2$. The parallelogram law then says $|z+i|^2+|z-i|^2=2|z|^2+2$, so we may factor $|z^2+1|^2-4|z|^2=|z^2+1|^2-2(2|z|^2+2)+4$ as $(|z+i|^2-2)(|z-i|^2-2)$. We thus see $z$ has norm at most $1$ and satisfies  
    \begin{equation}
        0\leq(|z+i|^2-2)(|z-i|^2-2).\label{eq: parallelogram law}\makelabel
    \end{equation}
    If the first factor vanishes, then $|z+i|^2=2$ and $|z-i|^2=2|z|^2+2-|z+i|^2=2|z|^2\leq 2$. Similarly, if the second factor vanishes, then $|z-i|^2=2$ and $|z+i|^2=2|z|^2+2-|z-i|^2=2|z|^2\leq 2$. That is, if $z$ lies on either disk's boundary, then it lies on the other disk and thus on the disks' intersection.

    Now say neither factor on the right of Equation $\eqref{eq: parallelogram law}$ vanishes, so both factors share a sign. Then either $z$ lies within both disks or neither. If $z$ lies outside both, we have $|z+i|^2>2$ and $|z-i|^2>2$. But then $2|z|^2+2=|z+i|^2+|z-i|^2>4$ and $z$ has norm more than $1$, contradicting our bond on $|z|$. Thus, for $z\in\Gamma$, then it must lie on the intersection of the two disks.

    Alternatively, if $z$ lies on both disks, then Equation $\eqref{eq: parallelogram law}$ is satisfied. Further, $|z+i|^2<2$ and $|z-i|^2<2$, so $2|z|^2+2=|z+i|^2+|z-i|^2<4$--which shows $|z|<1$. Thus, $\Gamma$ is the intersection of the disks of centers $\pm i$ and raddi $\sqrt{2}$.
    \end{proof}
    Lemma \ref{lem: propoerties of Cn} lets us later describe $\partial\Gamma$ as the solutions to $\frac{|z|}{|1+z^2|}=\frac{1}{2}$ with $|z|\leq 1$ or with a paramaterization.
See Figure $\ref{fig: Pn roots bound by Cn}$ to see both families of curves plotted for select values of $n$.

We now show $\partial\Gamma_n$ bounds the roots of $P_n(z)$ (see Figure $\ref{fig: Pn roots bound by Cn}$). 

\begin{theorem} \label{thm: Pn(z) has no roots in D_n}
    $P_n(z)$ has no roots in $\Gamma_n$ for all $n\geq 2$.
\end{theorem}

\begin{proof}
    Let $z\in \Gamma_n$. Writing $P_n(z)=(1+z)(1+z^2)^n-z^{n+1}R_n(z)$, dividing by $(1+z)(1+z^2)^n$, and taking moduli shows 
    $$\left|\frac{P_n(z)}{(1+z)(1+z^2)^n}\right|\geq1-\left|\frac{z^{n+1}R_n(z)}{(1+z)(1+z^2)^n}\right|\geq 1-K_n^n|zR_n(z)|,$$
    where the last inequality uses the definition of $\Gamma_n$. We show $K_n^n|zR_n(z)|<1$, thus giving the nonvanishing of $P_n(z)$ on $\Gamma_n$. The maximum modulus theorem bounds the values of $|zR_n(z)|$ within $\Gamma_n$ by any contour encasing it (avoiding singularities), so we consider the contour $|z|=1$ encasing $\Gamma_n$. As a polynomial of nonnegative coefficients, $zR_n(z)$ attains its maximum norm along the positive real axis, namely at $z=1$. Since this point lies in $\Gamma_n$, we see 
    $$K_n^n|zR_n(z)|\leq K_n^nR_n(1)=\left(
    \frac{n^2-1}{n^2}\right)^n<1$$
\end{proof}

\subsection{Uniform convergence to $\Gamma$}\label{subsec: convergence to Cn}
To prove the convergence of the roots of $P_n(z)$ to $\partial{\Gamma}_n$, we'll need a preliminary fact about $\sqrt[n]{R_n(z)}$.
\begin{lemma} 
    For all $0<t_0<1$, the sequence of functions $\sqrt[n]{R_n(z)}$ uniformly converges to $2$ on $|z|\leq t_0$.
\end{lemma}
\begin{proof}
    Fix $t_0\in(0,1)$ and consider $|z|\leq t_0$. We write $R_n(z)=\B{n}{0}A_n(z)$ where $A_n(z)=\sum_{k=0}^n a_{n,k}z^{k}$ and $a_{n,k}=\B{n}{k}/\B{n}{0}$. We show $A_n(z)$ uniformly converges to $\sum_{k\geq0}z^k=\frac{1}{1-z}$ for all $|z|\leq t_0$, so let $\epsilon>0$. We first fix some $M\in\ZZ^+$ such that $t_0^{M+1}<\frac{\epsilon}{2}(1-t_0)$ and consider $n$ already large enough for $M\leq n-1$. We now decompose $\frac{1}{1-z}-A_n(z)$ as
    \begin{equation}
        \sum_{k=0}^n(1-a_{n,k})z^k+\sum_{k=n+1}^{\infty}z^k=\sum_{k=0}^M(1-a_{n,k})z^k+\sum_{k=M+1}^n(1-a_{n,k})z^k+\sum_{k=n+1}^{\infty}z^k\makelabel\label{eq: sum decompositon}
    \end{equation}
    and bound the various sums individually. By maximizing the value of $(1-a_{n,k})$ in $k$ at $k=M$ and leveraging the triangle inequality,
    $$\left|\sum_{k=0}^M(1-a_{n,k})z^k\right|\leq (1-a_{n,M})\sum_{k=0}^M|z|^k\leq(1-a_{n,M})\left(\frac{1-t_0^{M+1}}{1-t_0}\right).$$
    We now show the value $(1-a_{n,M})$ may be made arbitrarily small. Particularly, we write $r_M=\lfloor n/2 \rfloor-\lfloor (n-M)/2 \rfloor$, so
   \begin{equation}
       a_{n,M}=\prod_{k=1}^{r_M}\frac{\lfloor n/2\rfloor-k+1}{\lfloor n/2\rfloor+k}=\prod_{k=1}^{r_M}\left(1-\frac{2k-1}{\lfloor n/2\rfloor+k}\right).\makelabel\label{eq (1-a) approx}
   \end{equation}
    For each $M$, this product has precisely $r_M$ terms. Then, for all $k\leq r_M$, each $\frac{2k-1}{\lfloor n/2\rfloor+j}$ term tends to zero as $n$ grows, so each factor of Equation $\eqref{eq (1-a) approx}$ tends to $1$. Thus the entire product in Equation $\eqref{eq (1-a) approx}$ tends to $1$, and so must $a_{n,M}$. Since the term $(1-a_{n,M})$ may then be made arbitrarily small for fixed $M$, pick some $N$ such that for all $n\geq N$, the inequality $(1-a_{n,M})<\frac{\epsilon}{2}\left(\frac{1-t_0}{1-t_0^{M+1}}\right)$ holds. Then we see
    \begin{equation}
        \left|\sum_{k=0}^n(1-a_{n,k})z^k\right|\leq (1-a_{n,M})\left(\frac{1-t_0^{M+1}}{1-t_0}\right)< \frac{\epsilon}{2}.\makelabel\label{eq: first sum bound}
    \end{equation}
The remaining two sums of Equation $\eqref{eq: sum decompositon}$, viewed as one sum in $k\geq M+1$, has all coefficients with norm at most $1$, so
    \begin{equation}  
    \left|\sum_{k=M+1}^n(1-a_{n,k})z^k+\sum_{k=n+1}^{\infty}z^k\right|<\left|\sum_{k=M+1}^\infty z^k\right|=\frac{t_0^{M+1}}{(1-t_0)}<\frac{\epsilon }{2}.\makelabel\label{eq: second sum bound}
    \end{equation}
    Combining Equations $\eqref{eq: first sum bound}$ and $\eqref{eq: second sum bound}$ in Equation $\eqref{eq: sum decompositon}$, we see that for all $n\geq N$, 
    \begin{equation}
        \left|\sum_{k=0}^\infty z^k -A_n(z)\right|<\frac{\epsilon}{2}+\frac{\epsilon}{2}=\epsilon.
    \end{equation}
    Then for all $\epsilon>0$, there exists some $N$ such that for all $n\geq N$ and $|z|<t_0$,
    \begin{align*}
        A_n(z)\leq \frac{1}{1-z}+\epsilon\leq \frac{1}{1-t_0}+\epsilon\\
        A_n(z)\geq \frac{1}{1-z}-\epsilon\geq \frac{1}{1-t_0}-\epsilon
    \end{align*}    
    By setting $c=\frac{1}{1-t_0}-\epsilon$ and $C=\frac{1}{1-t_0}+\epsilon$ we see $c\B{n}{0}\leq R_n(z)
    \leq C\B{n}{0}$ and 
    \begin{equation}
        c^{\frac{1}{n}}\B{n}{0}^{\frac{1}{n}}\leq \sqrt[n]{R_n(z)}
    \leq C^{\frac{1}{n}}\B{n}{0}^{\frac{1}{n}}.
    \end{equation}
    Since $n$th roots of real constants tend to $1$, the squeeze theorem limits $\sqrt[n]{R_n(z)}$ to the limit of $\B{n}{0}^{\frac{1}{n}}$, which we show is $2$.
    
    We present a more combinatorial argument, but could alternatively appeal to Stirling's approximation. First, the SDT of shape $(n,n)$, counted by $\B{n}{0}$ is a strict subset of all 2 row $n$ domino SDT, counted by $2^n$. This shows $\B{n}{0}\leq2^n$, and we deduce $\B{n}{0}^{\frac{1}{n}}<2$. On the other hand, the number of SDT of shape $(n+k,n-k)$ is maximized at $k=0$ by Theorem \ref{thm: 2 col SDT satisfy pascalian recursion}, so
    
    $$2^n=\sum_{k=0}^n\B{n}{k}<(n+1)\B{n}{0},$$
    which shows $\B{n}{0}^{\frac{1}{n}}>\frac{2}{\sqrt[n]{1+n}}$ with limit $2$. By the squeeze theorem, $\B{n}{0}^{\frac{1}{n}}$, and $\sqrt[n]{R_n(z)}$ for $|z|<t_0$, converges to $2$.
\end{proof}
The condition that $t_0<1$ cannot be done away with, since the sequence $R_n(-1)^{\frac{1}{n}}$ has no limit: Particularly, $R_{2m+1}(-1)=0$ for all $m$ while $\lim_{n\to\infty}R_{2m}(-1)^{1/2m}=\lim_{n\to\infty}\B{2m}{0}^{1/2m}=2$. In this sense, our convergence result is optimal.
\begin{theorem}\label{thm: PP roots converge to C_infty}
    The roots of $P_n(z)$ converge uniformly to $\partial\Gamma$ and fill the curve densely.
\end{theorem}
\begin{proof}
    Theorem $\ref{thm: PP roots annuli}$ lets us limit our view to the complex unit disk. Following in the proof of Theorem $\ref{thm: Pn(z) has no roots in D_n}$ above, the solutions of $P_n(z)$ satisfy

    $$1=\frac{z^{n+1}R_n(z)}{(1+z)(1+z^2)^n}.$$
    Since no Pascalian polynomial has positive real zeros, we take $n$th roots with a branch cut along the positive real axis and write the desired zeros as solutions to 
    $$1=\left(\frac{z\sqrt[n]{R_n(z)}}{(1+z^2)}\sqrt[n]{\frac{z}{1+z}}\right)^n.$$
    Since we're examining the roots of $P_n(z)$, we exclude a neighborhood around $z=0$ from our view, and let $\sqrt[n]{\frac{z}{z+1}}$ converge uniformly to $1$. Similarly, we use the uniform convergence of $\sqrt[n]{R_n(z)}$ to $2$ to see the roots of $P_n(z)$ are asymptotically the solutions to
    $$1=\left(\frac{2z}{1+z^2}\right)^{n}.$$
    The map $w=2z/(1+z^2)$ maps all points of $\partial\Gamma$ to the complex unit circle, and, given $|w|=1$, we see $\frac{|z|}{|1+z^2|}=\frac{1}{2}$. Then solving for $z$ in $w$ shows
    \begin{equation}
        wz^2-2z+w=0\makelabel\label{eq: quadratic}
    \end{equation}
    By appealing to Vieta's relations, we see the solutions multiply to $1$, so at least one solution obeys $|z|\leq1$ and thus lies on $\partial\Gamma$. If both roots have norm $1$, they're either both $1$ or $-1$, the only $z\in\partial\Gamma_n$ with $|z|=1$. Only the former yields a solution to equation $\eqref{eq: quadratic}$, in which case $z=1$ is a double root. This shows any $w$ on the complex unit circle has a unique preimage on $\Gamma$, so $w$ is a bijection of $\Gamma$ to the complex unit circle. As a consequence, there exist points $z_m\in\partial\Gamma$ for $1\leq m\leq n$ such that
    
    $$\frac{2z_m}{1+z_m^2}=e^{2\pi i m/n}.$$
    By the asymptotic density of the $n$th roots of unity on the unit circle, $z_m$ are asymptotically dense on the curve $\Gamma$. As we've shown their existence on $\Gamma$ and that they're asymptotically the roots of $P_n(z)$, we're done with our statement on $P_n(z)$. 
\end{proof}
Figure $\ref{fig: zm plotted n=7,21,35}$ plots the roots of $P_n(z)$ and their approximations $z_m$ for select values of $n$. 

\begin{figure}[htbp]
  \centering  \hspace*{.7in}
  \begin{minipage}[t]{0.22\linewidth}
    \centering
    \includegraphics[width=\linewidth]{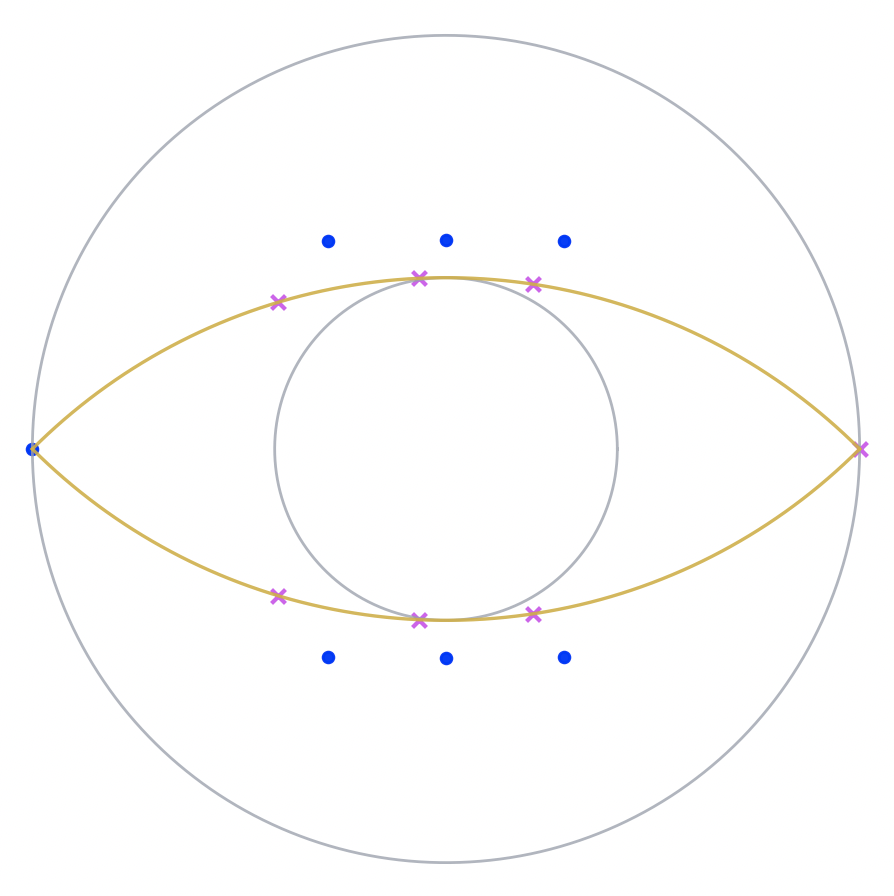}
  \end{minipage}\hfill
  \begin{minipage}[t]{0.22\linewidth}
    \centering
    \includegraphics[width=\linewidth]{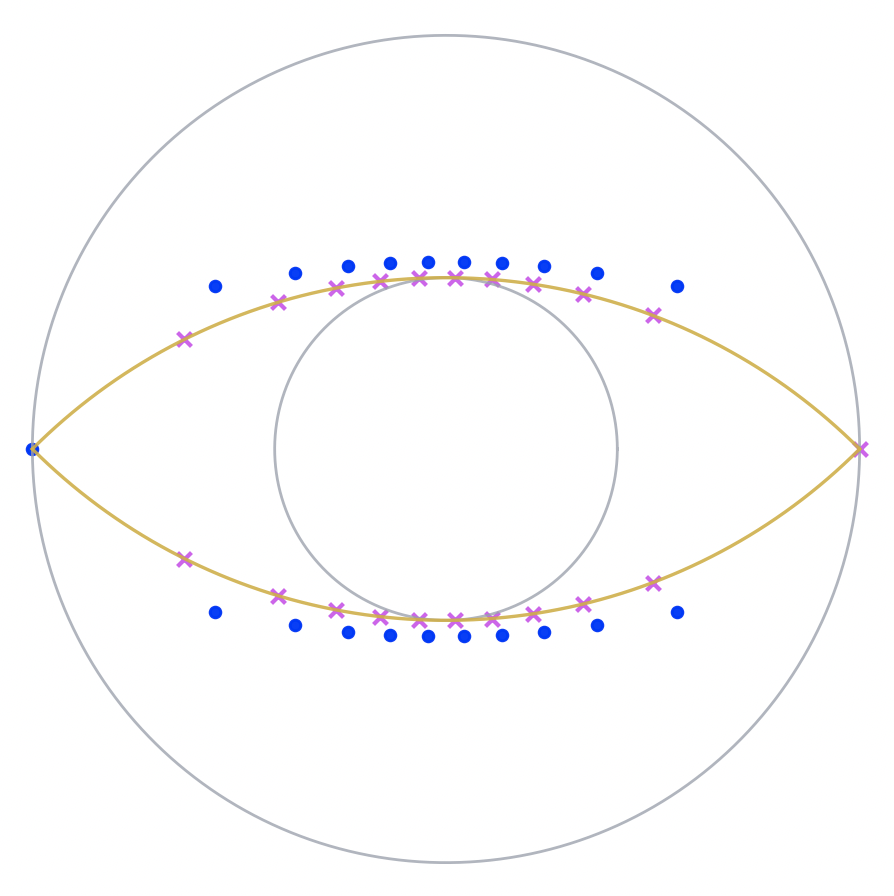}
  \end{minipage}\hfill
  \begin{minipage}[t]{0.22\linewidth}
    \centering
    \includegraphics[width=\linewidth]{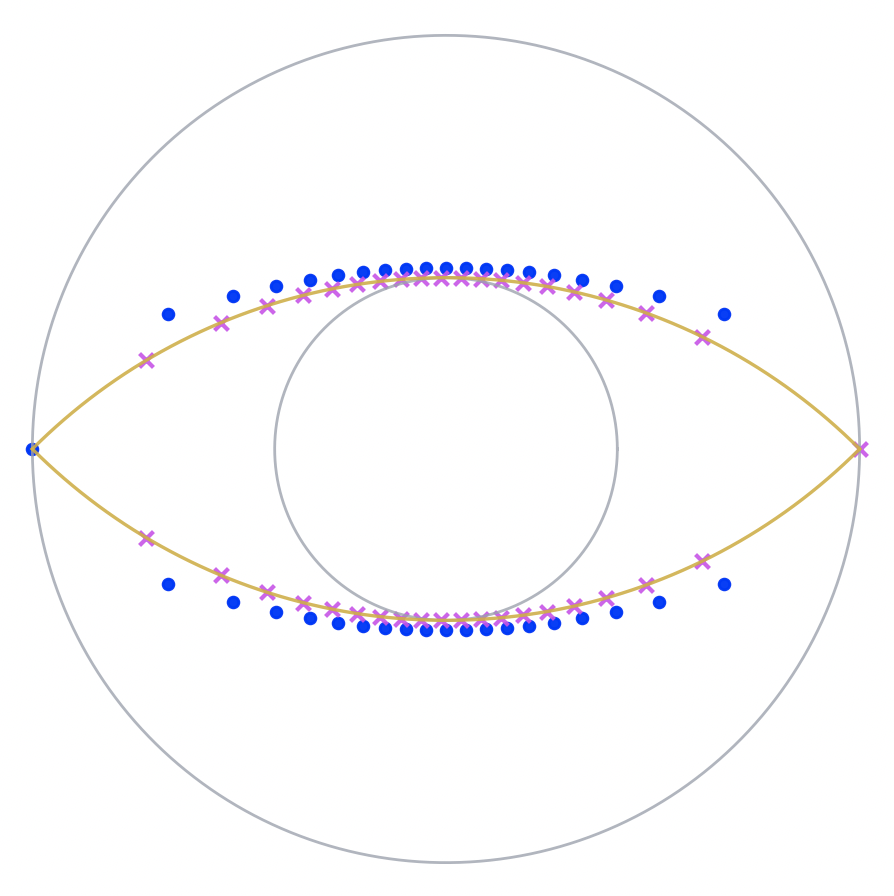}
  \end{minipage}  \hspace*{.7in}
  \caption{The roots of $P_n(z)$ (blue) with their  approximations $z_m$ (pink) on $\partial\Gamma$ (gold) for $n=7,21,35$. The roots of $P_n(z)$ converge to the points $z_m$ and are asymptotically dense in $\partial\Gamma$.}
  \label{fig: zm plotted n=7,21,35}
\end{figure}

\section{Irreducibility and Galois Groups}\label{sec: irreducibility and galois groups}
We use this section to briefly discuss the irreducibility and Galois groups of the Pascalian polynomials. While we derive the prime factorization of $P_n(z)$ for odd $n$ using work on the irreducibility of truncated binomial polynomials \cite{DubickasSiurys2017,FilasetaKumchevPasechnik2007,KhandujaKhassaLaishram2011,laishram-yadav2024}, we have only conjectures for the irreducibility of the even indexed Pascalian polynomials. We similarly make certain deductions of $P_n(z)$'s Galois group for odd $n$ and form a conjecture on $P_n(z)$'s Galois group in general. 

\subsection{Irreducibility}\label{subsec: irreducibility}
 While the general problem of the irreducibility of truncated binomial polynomials remains open, the irreducibility of the center truncation $q_n(z)$ has been given several times (see  \cite{KhandujaKhassaLaishram2011,laishram-yadav2024}, for instance). To deduce the factorization of $P_n(z)$ into irreducibles for odd $n$, we first prove the pleasant lemma below, which could not be verified as known. We state the argument in general for its applications outside our context.

\begin{lemma}\label{lem: q(z^2) criteria}
    Let $R$ be a domain and consider an irreducible $q(z)\in R[z]$. If $q(z^2)$ is reducible in $R[z]$, it must have square constant and leading coefficients. 
\end{lemma}

\begin{proof}
    Let $q(z)$ be irreducible, and say $q(z^2)$ factors into irreducibles $g_1(z)g_2(z)\dots g_n(z)$. No $g_i(z)$ may be even, since then $q(z)$ would be reducible with polynomial divisor $g_i(\sqrt{z})$. Setting $z\mapsto -z$ in $q(z^2)$'s factorization shows
    $$g_1(-z)g_2(-z)\dots g_n(-z)=q\left((-z)^2\right)=q\left(z^2\right)=g_1(z)g_2(z)\dots g_n(z).$$
    We deduce, $n$ is even and the factors of $q(z^2)$ may be partitioned into $\{g_i(z),g_i(-z)\}$ pairs. Letting $h(z)$ be the product of $n/2$ of the $g_i(z)$, each from a distinct pair, lets us write $q(z^2)=h(z)h(-z)$. Thus, the leading and constant terms of $q(z^2)$ are the squares of the leading and constant terms of $h(z)$ respectively.
\end{proof}
We deduce the irreducibility of $q_n(z^2)$ by showing $\B{n}{n}=\binom{n}{\lfloor n/2\rfloor}$, the leading coefficient of $P_n(z)$, is not a square for odd $n\geq 3$.
\begin{theorem}\label{thm: P_n(z) prime factorization for odd n}
    For odd $n$, $P_n(z)$ factors into irreducibles over $\QQ$ as $(1+z)q_n(z^2)$.
\end{theorem}

\begin{proof}
    Consider odd $n=2m+1$. $q_n(z)$ is irreducible with leading coefficient $\binom{2m+1}{m}=\frac{(2m+1)!}{m!(m+1)!}$. Bertrand's postulate promises a prime $p$ such that $m+1<p<2(m+1)$. This prime divides the numerator of $\frac{(2m+1)!}{m!(m+1)!}$ precisely once and the denominator not at all, so $\binom{2m+1}{m}$ is never a square and $q_n(z^2)$ must be irreducible.
\end{proof}
We alternatively have the following conjecture, affirmed by numerical computations for low $n$.
\begin{conjecture}
    $P_n(z)$ is irreducible over $\QQ$ for all even $n$.
\end{conjecture}
\subsection{Galois Groups}\label{subsec: galois groups}
For $n=2m+1$, the Galois group of $P_n(z)$ over $\QQ$ is precisely the Galois group of $q_n(z^2)$. It is well known (see \cite{Odoni}) that the Galois group of the polynomial $f(z^2)$ embeds inside the wreath product $\ZZ_2\wr \operatorname{GAL}f(z)$, so our attention is again brought back to the truncated binomial expansions. It has been shown that $q_n(z)$ has full Galois group for all $n$ sufficiently large, but it is still only conjectured that $q_n(z)$ (and all nontrivial truncated binomial expansions) has full Galois group. Most recently, Laishram and Yadav \cite{laishram-yadav2024} showed that a wide family of truncated binomial polynomials all have full Galois group under the assumption of the exploit $abc$ conjecture.

\begin{conjecture}\label{conj: gal groups}
    $P_{n}(z)$ has full Galois group for even $n$. For $n=2m+1$, the Galois group of $P_n(z)$ is the $m$th hyperoctahedral group.
\end{conjecture}

For $P_{2m+1}(z)$ to have full Galois group as conjectured, it is necessary but not sufficient that $q_n(z)$ have full Galois group.

\section{Directions for Further Study}\label{sec: directions for further study}
In Section $\ref{sec: roots of PP}$, we used $P_n(z)$'s recursion to show $P_n(z)$ shares no (nontrivial) roots with $P_{n-2}(z)$. Can $P_n(z)$ and $P_{m}(z)$ share nontrivial roots for general $n\neq m$?  Can the extended recursion of Theorem $\ref{thm: extended Pn recursion, general k}$ be used to examine the common roots of $P_n(z)$ and $P_{n-k}(z)$ as with $k=2$?

In Section $\ref{sec: limit curve of roots}$, we classified the root asymptotics of a family of polynomials made by sorting the coefficients of a given family of polynomials. The roots of polynomials whose coefficients are the sorted Eulerian numbers, Naryana numbers, and type $B$ Narayna numbers all seem to have similar limit curves to $\partial\Gamma$. Do analogous results hold for other families? Does a more general setting exist for families of polynomials whose coefficients obey a normal distribution? What about other distributions?

Section $\ref{sec: irreducibility and galois groups}$ made limited progress on $P_n(z)$'s algebra. A natural extension would be an examination of $P_n(z)$'s irreducibility for even $n$ and a classification of $P_n(z)$'s Galois groups.

\section{Acknowledgments}
I'd like to thank Dr. Allan Berele for his extensive support and mathematical guidance. He has been an attentive mentor, a dedicated collaborator, and a good friend over the course of my mathematical career.

\end{document}